\documentclass[twoside]{article} 
\usepackage{graphicx}
\date{} 
\title{The numerical evaluation of the Riesz function}
\author{\sc R. B.\ Paris \\
{\em Division of Computing and Mathematics,} \\
{\em Abertay University, Dundee DD1 1HG, UK}}
\begin{document}
\def\f#1#2{\mbox{${\textstyle \frac{#1}{#2}}$}}
\def\dfrac#1#2{\displaystyle{\frac{#1}{#2}}}
\def\boldal{\mbox{\boldmath $\alpha$}}
\newcommand{\bee}{\begin{equation}}
\newcommand{\ee}{\end{equation}}
\newcommand{\sa}{\sigma}
\newcommand{\ka}{\kappa}
\newcommand{\al}{\alpha}
\newcommand{\la}{\lambda}
\newcommand{\ga}{\gamma}
\newcommand{\eps}{\epsilon}
\newcommand{\om}{\omega}
\newcommand{\fr}{\frac{1}{2}}
\newcommand{\fs}{\f{1}{2}}
\newcommand{\g}{\Gamma}
\newcommand{\br}{\biggr}
\newcommand{\bl}{\biggl}
\newcommand{\ra}{\rightarrow}
\newcommand{\gtwid}{\raisebox{-.8ex}{\mbox{$\stackrel{\textstyle >}{\sim}$}}}
\newcommand{\ltwid}{\raisebox{-.8ex}{\mbox{$\stackrel{\textstyle <}{\sim}$}}}
\renewcommand{\topfraction}{0.9}
\renewcommand{\bottomfraction}{0.9}
\renewcommand{\textfraction}{0.05}
\newcommand{\mcol}{\multicolumn}
\date{}
\maketitle
\pagestyle{myheadings}
\markboth{\hfill \sc R. B.\ Paris  \hfill}
{\hfill \sc Asymptotics of $F_{n,\sa}(x;\mu)$\hfill}
\begin{abstract}
The behaviour of the generalised Riesz function defined by
\[S_{m,p}(x)=\sum_{k=0}^\infty \frac{(-)^{k-1}x^k}{k! \zeta(mk+p)}\qquad (m\geq 1,\ p\geq 1)\]
is considered for large positive values of $x$. A numerical scheme is given to compute this function which enables the visualisation of its asymptotic form. The two cases $m=2$, $p=1$ and $m=p=2$ (introduced respectively by Hardy and Littlewood in 1918 and Riesz in 1915) are examined in detail. It is found on numerical evidence that these functions appear to exhibit the $x^{-1/4}$ and $x^{-3/4}$ decay, superimposed on an oscillatory structure, required for the truth of the Riemann hypothesis. 
\vspace{0.3cm}

\vspace{0.1cm}
 
\end{abstract}

\vspace{0.3cm}

\noindent $\,$\hrulefill $\,$

\vspace{0.3cm}

\begin{center}
{\bf 1.\ Introduction}
\end{center}
\setcounter{section}{1}
\setcounter{equation}{0}
\renewcommand{\theequation}{\arabic{section}.\arabic{equation}}
The generalised Riesz function is defined by the sum
\bee\label{e11}
S_{m,p}(x)=\sum_{k=0}^\infty \frac{(-)^{k-1} x^k}{k! \zeta(mk+p)}\qquad (x\geq0),
\ee
where $m\geq 1$, $p\geq 1$ and $\zeta(s)$ is the Riemann zeta function. The original function considered by Riesz \cite{R} took the form
\[\sum_{k=1}^\infty\frac{(-)^{k-1} x^k}{(k-1)! \zeta(2k)},\]
which corresponds to a case of (\ref{e11}) with $m=p=2$ since it is easily seen to equal $-x S_{2,2}(x)$. A similar function corresponding to $m=2$, $p=1$ was discussed in the famous memoir by Hardy and Littlewood \cite{HL}. The interest in both these cases of (\ref{e11}) results from the fact that a necessary and sufficient condition for the truth of the Riemann hypothesis is that \cite[p.~382]{T}
\bee\label{e12}
S_{2,1}(x)=O(x^{-1/4+\epsilon}),\qquad S_{2,2}(x)=O(x^{-3/4+\epsilon})
\ee
as $x\to+\infty$, where $\epsilon$ is an arbitrarily small positive quantity. The results in (\ref{e12}) are superficially attractive as they are derived from sums containing only values of $\zeta(s)$ at positive integer values of $s$.

The sum in (\ref{e11}) can also be viewed as an example of a perturbation of the exponential series for $e^{-x}$ in the form $\sum_{k\geq 0} a_k (-x)^k/k!$, where $a_k$ are coefficients that possess the property $a_k\to1$ as $k\to\infty$; in the case of the Riesz function we have $a_k=1/\zeta(mk+p)$. The growth of this series for large (complex) $x$ is found to depend sensitively on the decay of the perturbing coefficients $a_k$. A discussion of this problem, together with several examples, is given in \cite{P}.

This paper is partly based on the earlier report by the author in \cite{PRiesz}. We present a computational scheme
for the numerical evaluation of $S_{m,p}(x)$ for large positive values of $x$. In particular, we concentrate on the Hardy-Littlewood case of $m=2$, $p=1$ and also on the Riesz case $m=p=2$ and determine numerically their large-$x$ behaviour. Based on our numerical results, we conclude that these cases are characterised by a damped oscillatory structure for sufficiently large positive $x$ with an amplitude that corresponds to the estimates in (\ref{e12}).

\vspace{0.6cm}

\begin{center}
{\bf 2.\ A computational scheme}
\end{center}
\setcounter{section}{2}
\setcounter{equation}{0}
\renewcommand{\theequation}{\arabic{section}.\arabic{equation}}
We use the result \cite[p.~3]{T}
\[\frac{1}{\zeta(s)}=\sum_{n=1}^\infty\frac{\mu(n)}{n^s}\qquad (\Re (s)>1),\]
where $\mu(n)$ is the M\"obius function defined by $\mu(n)=(-)^r$ if $n$ has $r$ distinct primes (with $\mu(1)=1$) and $\mu(n)=0$ otherwise. Then we obtain the expansion for the case $m=2$, $p\geq1$
\begin{eqnarray}
S_{2,p}(x)\!\!&=&\!\! \sum_{k=0}^\infty \frac{(-)^{k-1}x^k}{k! \zeta(2k+p)}\nonumber\\
&=&\!\!\sum_{n=1}^\infty \frac{\mu(n)}{n^p} (1-e^{-X_n})-\frac{1}{\zeta(p)},\qquad X_n:=\frac{x}{n^2}.\label{e21}
\end{eqnarray}
The factor $1-e^{-X_n}$ has a ``cut-off'' when $n=N=\lceil x^{1/2}\rceil$; for $n>N$, the decay of the terms in this series is slow with late terms eventually behaving like $\mu(n)/n^{p+2}$.

This slow decay in the tail can be accelerated as follows. We write
\bee\label{e22}
S_{2,p}(x)=\sum_{n=1}^{N-1}\frac{\mu(n)}{n^p} (1-e^{-X_n})+T_N(p;x)-\frac{1}{\zeta(p)},\qquad T_N(p;x):=\sum_{n=N}^\infty \frac{\mu(n)}{n^p} (1-e^{-X_n}),
\ee
where $N$ is chosen as above. In the tail $T_N(p;x)$ (where $X_n\leq 1$) we put
\[1-e^{-X_n}=X_n\bl(1-\frac{X_n}{2!}+\frac{X_n^2}{3!}-\cdots\br)=X_n \,{}_1F_1(1;2;-X_n)\]
\[=X_n e^{-X_n} {}_1F_1(1;2;X_n)=X_n f_1(X_n),\]
where 
\bee
f_k(z):=e^{-z}{}_1F_1(1;k+1;z)\qquad (k\geq0).
\ee
Here ${}_1F_1$ denotes the confluent hypergeometric function and Kummer's transformation \cite[p.~325]{DLMF} has been used to change the argument from $-X_n$ to $X_n$.

From (\ref{a2}) in the appendix the following lemma is established:
\medskip

\noindent{\bf Lemma 1.}\ \ The function $zf_1(z)$ satisfies the recursion formula
\bee\label{e23}
zf_1(z)=1-\frac{1}{e_k(z)}+\frac{z^k f_k(z)}{k!\,e_k(z)}\qquad (k\geq 1),
\ee
where $e_k(z)$ denotes the sum of the first $k$ terms of the exponential series
\bee\label{e24}
e_k(z):=\sum_{r=0}^{k-1}\frac{z^r}{r!}.
\ee
\medskip

From the definition (\ref{e24}) it follows that
\bee\label{e25}
\frac{1}{e_k(z)}=e_k(-z)+\frac{z^k g_k(z)}{k!\,e_k(z)},
\ee
where $g_k(z)$ is a polynomial of degree $k-2$. Some routine algebra shows that
\[g_2(z)=2,\quad g_3(z)=-\frac{3}{2}z,\quad g_4(z)=2+\frac{2}{3}z^2,\]
\[g_5(z)=-\frac{5}{3}z-\frac{5}{24}z^3,\quad g_6(z)=2+\frac{3}{4}z^2+\frac{1}{20}z^4, \ldots\, .\]

Then, using (\ref{e23})--(\ref{e25}), the tail of the series can be written in the form
\bee\label{e26}
T_N(p;x)=T_N^{(1)}(p;x)+T_N^{(2)}(p;x),
\ee
where
\bee\label{e27}
T_N^{(1)}(p;x)=
\frac{1}{k!}\sum_{n=N}^\infty\frac{\mu(n)}{n^p}\,\frac{X_n^k\,\Delta_k(X_n)}{e_k(X_n)},\qquad \Delta_k(X_n):=f_k(X_n)-g_k(X_n)
\ee
and
\bee\label{e28}
T_N^{(2)}(p;x)=\sum_{n=N}^\infty\frac{\mu(n)}{n^p}\sum_{r=1}^{k-1}\frac{(-)^{r-1}X_n^r}{r!}=\sum_{r=1}^{k-1}\frac{(-X_N)^r}{r!}\,\lambda_r
\ee
with
\[X_N:=\frac{x}{N^2},\qquad\lambda_r:=N^{2r}\bl\{\sum_{r=1}^{N-1}\frac{1}{n^{2r+p}}-\frac{1}{\zeta(2r+p)}\br\}.\]
Since $X_n\to0$ as $n\to\infty$, we see that the decay of the late terms in $T_N^{(1)}(x)$ is now controlled by $\mu(n)/n^{2k+p}$, which for $k\geq 2$ represents a modest improvement in the rate of convergence of the series.

Collecting together the results in (\ref{e22}), (\ref{e26}) -- (\ref{e28}), we have
\bee\label{e29}
S_{2,p}(x)=\sum_{n=1}^{N-1}\frac{\mu(n)}{n^p} (1-e^{-X_n})-\frac{1}{\zeta(p)}+T_N^{(1)}(p;x)+T_N^{(2)}(p;x)
\ee
for $p\geq 1$.

\vspace{0.6cm}

\begin{center}
{\bf 3.\ Numerical results for $m=2$ and $p=1, 2$}
\end{center}
\setcounter{section}{3}
\setcounter{equation}{0}
\renewcommand{\theequation}{\arabic{section}.\arabic{equation}}
We have employed the scheme (\ref{e29}) with $k=6$ to compute the Hardy-Littlewood case $m=2$ $p=1$ for $x\geq0$ up to $x=10^8$. The results are shown in the sequence of plots in Fig.~1. It is found that $S(x)\equiv S_{2,1}(x)$ decreases once $x\,\gtwid\,3$ down to values of the order $10^{-6}$, whereupon the graph commences to oscillate about the zero line. The oscillations appear to be regular and have a decreasing amplitude.
\begin{figure}[th]
	\begin{center}	{\tiny($a$)}\includegraphics[width=0.4\textwidth]{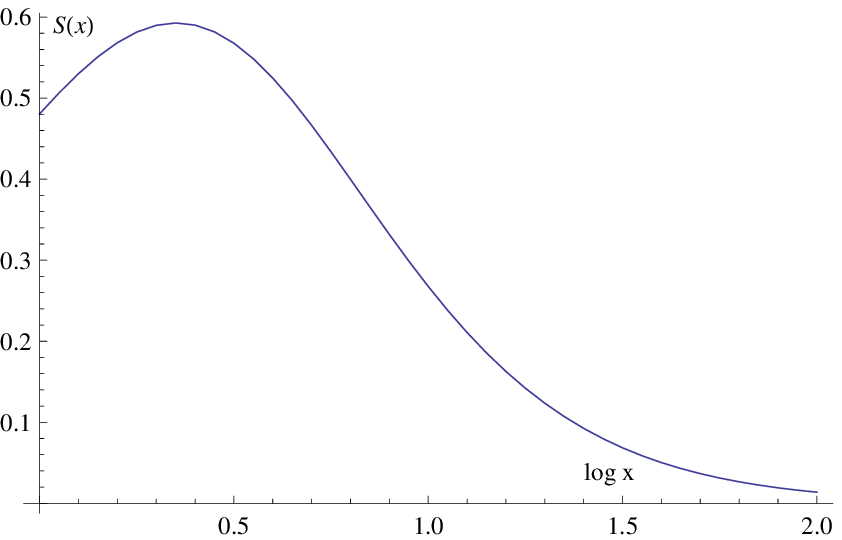}
	\qquad
	{\tiny($b$)}\includegraphics[width=0.4\textwidth]{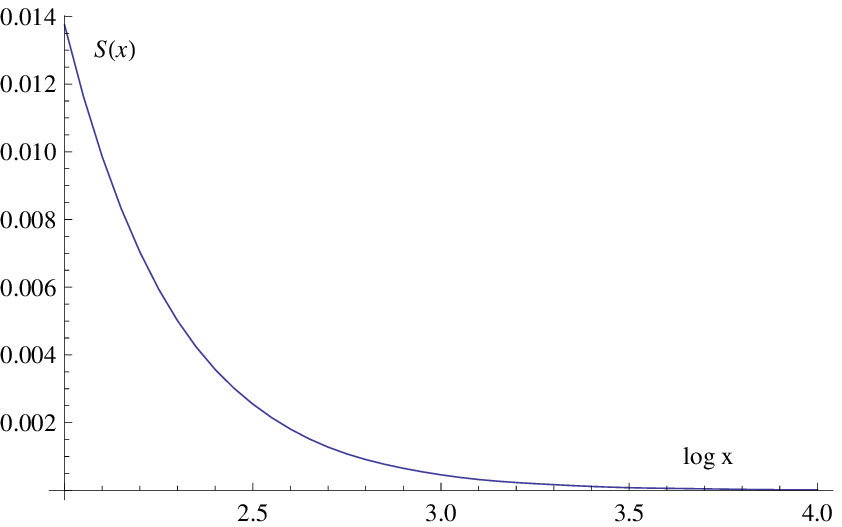} 
	\bigskip
	
		{\tiny($c$)}\includegraphics[width=0.4\textwidth]{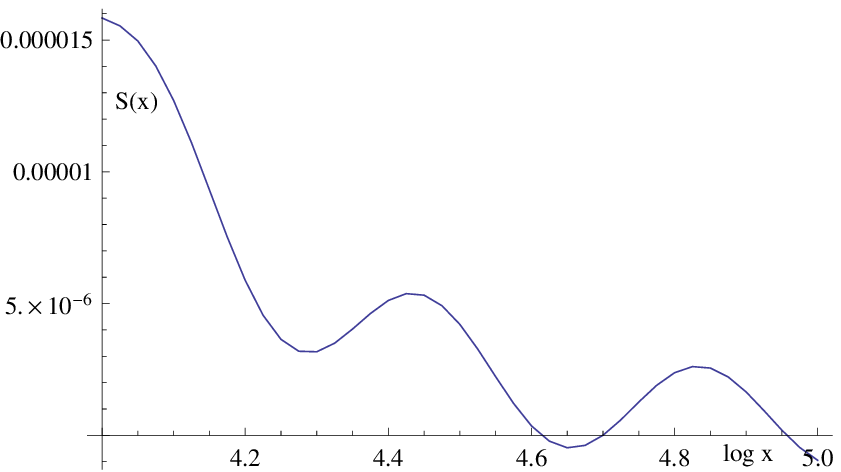}
	\qquad
	{\tiny($d$)}\includegraphics[width=0.4\textwidth]{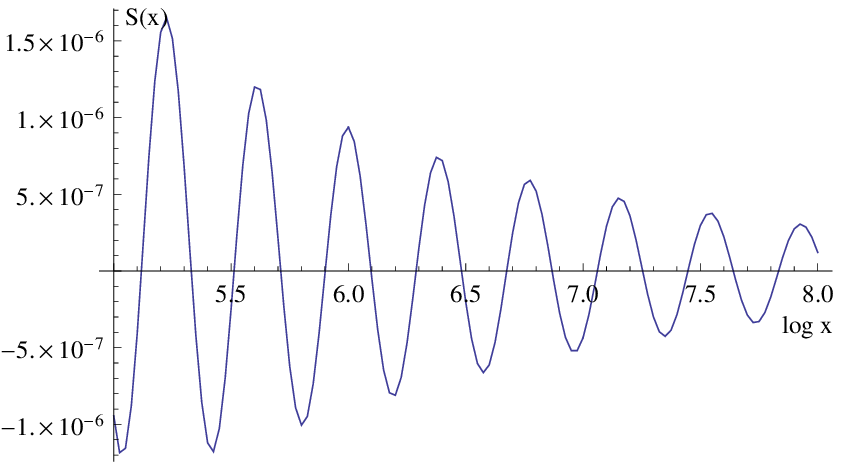} \\
	
\caption{\small{The graph of $S(x)\equiv S_{2,1}(x)$ against $\log_{10}x$ for different ranges of $x$.}}
\end{center}
\end{figure}

The successive maxima\footnote{We commence the enumeration of the maxima and minima from the point where the graph of $S(x)$ first becomes negative.} ($x_k^+$) and minima ($x_k^-$) in the oscillatory region are determined by calculating the zeros of the derivative $S_{2,1}'(x)$, where
\[S_{2,1}'(x)=\sum_{n=1}^\infty \frac{\mu(n)}{n^3} \,e^{-X_n}=\frac{1}{\zeta(3)}-\sum_{n=1}^\infty \frac{\mu(n)}{n^3} \,(1-e^{-X_n}),\]
by using (\ref{e29}) with $p=3$. 
The results of these calculations together with the corresponding values of $S(x_k^\pm)$ are shown in Table 1. Plots of $\log\,S_{2,1}(x_k^\pm)$ against $\log\,x$ are shown in Fig.~2, where it is seen that they reveal a linear variation to a good approximation. The dashed lines in these figures have slope equal to $-0.25$, thereby numerically confirming the estimate in (\ref{e12}).
\begin{table}[h]
\caption{\footnotesize{Values of the maxima and minima $\log_{10} x_k^\pm$ and the corresponding values of $S_{2,1}(x_k^\pm)$.}}
\begin{center}
\begin{tabular}{|c|ll|ll|}
\hline
&&&&\\[-0.3cm]
\mcol{1}{|c|}{$k$} & \mcol{1}{c}{$\log_{10} x_k^+$} & \mcol{1}{c|}{$S_{2,1}(x_k^+)$} &\mcol{1}{c}{$\log_{10} x_k^-$} & \mcol{1}{c|}{$S_{2,1}(x_k^-)$} \\
[.1cm]\hline
&&&&\\[-0.25cm]
1 & 4.83284 & $2.62044\times 10^{-6}$ & 4.65573 & $-4.78520\times 10^{-7}$ \\
2 & 5.22278 & $1.65298\times 10^{-6}$ & 5.03476 & $-1.20507\times 10^{-6}$ \\
3 & 5.61033 & $1.21568\times 10^{-6}$ & 5.41918 & $-1.18328\times 10^{-6}$ \\
4 & 5.99669 & $9.38573\times 10^{-7}$ & 5.80406 & $-1.00735\times 10^{-6}$ \\
5 & 6.38315 & $7.47310\times 10^{-7}$ & 6.19039 & $-8.20165\times 10^{-7}$ \\
6 & 6.76905 & $5.93349\times 10^{-7}$ & 6.57596 & $-6.62603\times 10^{-7}$ \\
7 & 7.15545 & $4.76797\times 10^{-7}$ & 6.96249 & $-5.30007\times 10^{-7}$ \\
8 & 7.54124 & $3.80007\times 10^{-7}$ & 7.34817 & $-4.26057\times 10^{-7}$ \\
\hline
\end{tabular}
\end{center}
\end{table}
\begin{figure}[th]
	\begin{center}	{\tiny($a$)}\includegraphics[width=0.4\textwidth]{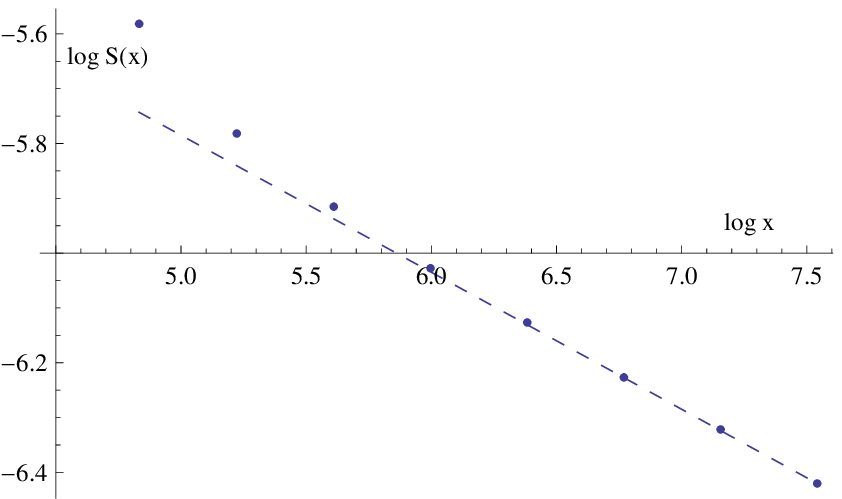}
	\qquad
	{\tiny($b$)}\includegraphics[width=0.4\textwidth]{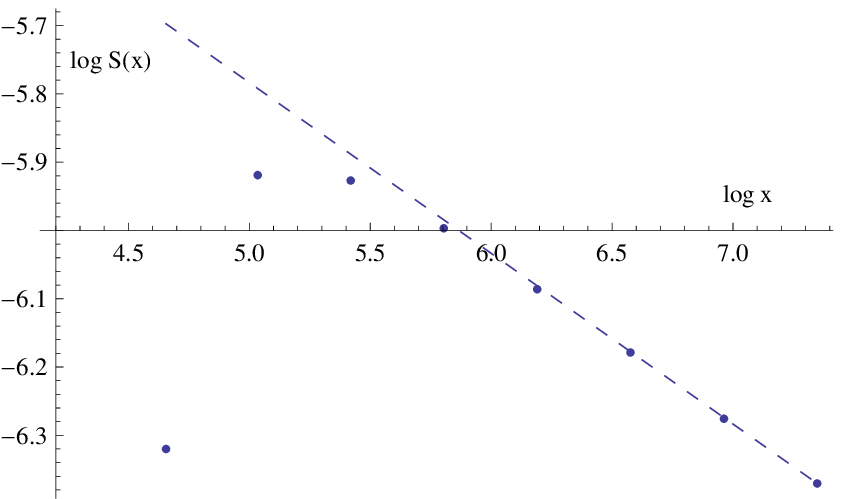} 
\caption{\small{Plots of (a) $\log_{10} S(x_k^+)$ and (b) $\log_{10} S(x_k^-)$ against $\log_{10}x$ (where $S(x)\equiv S_{2,1}(x)$). The dashed lines have slope $-0.25$.}}
\end{center}
\end{figure}
\begin{figure}[th]
	\begin{center}	{\tiny($a$)}\includegraphics[width=0.4\textwidth]{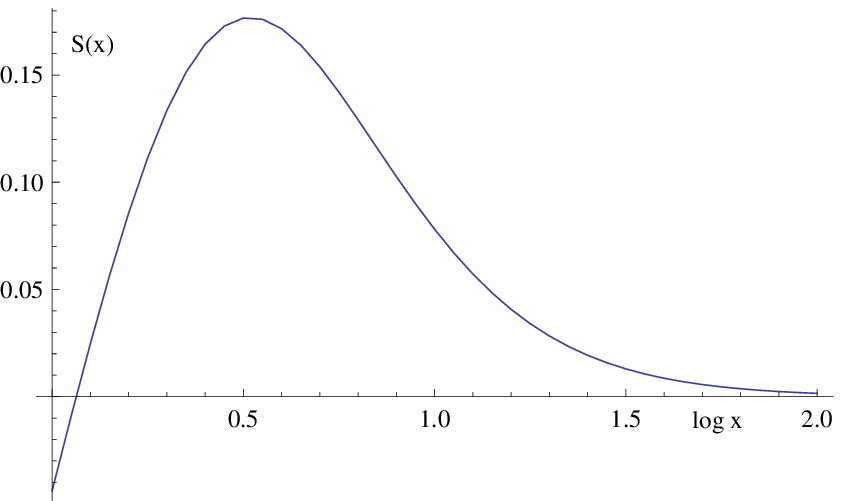}
	\qquad
	{\tiny($b$)}\includegraphics[width=0.4\textwidth]{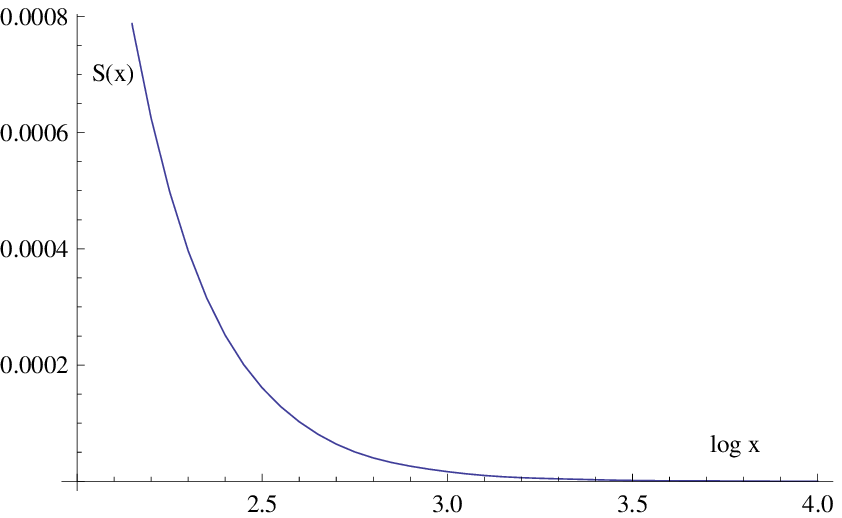} 
	\bigskip
	
		{\tiny($c$)}\includegraphics[width=0.4\textwidth]{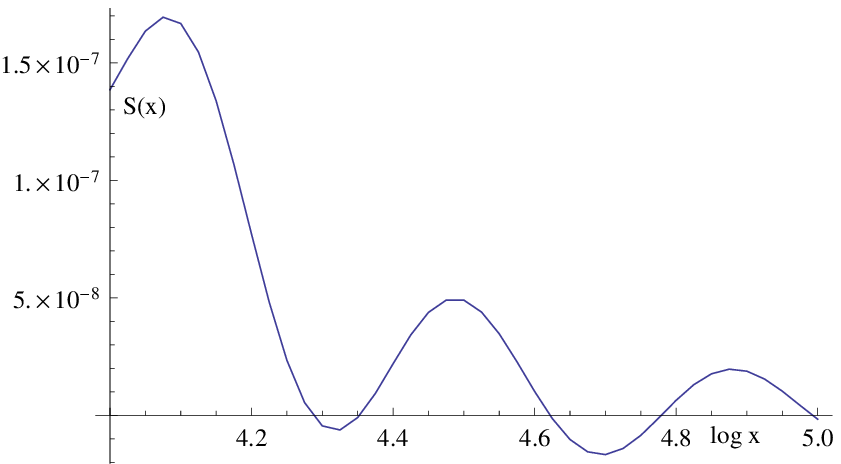}
	\qquad
	{\tiny($d$)}\includegraphics[width=0.4\textwidth]{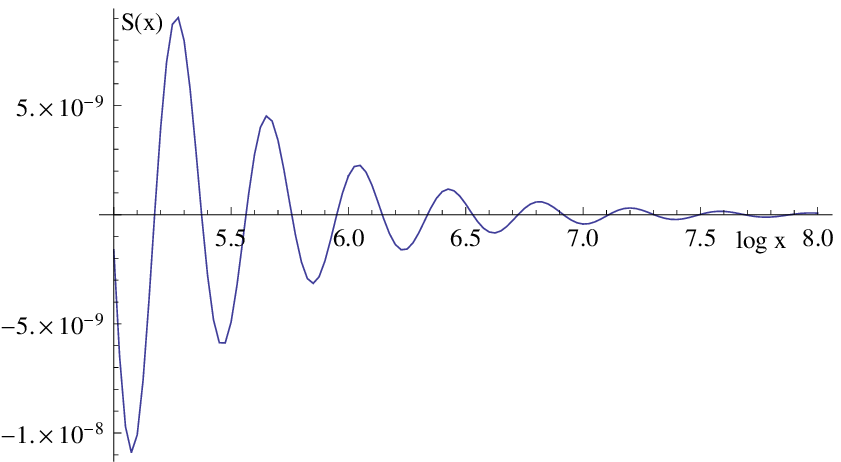} \\
	
\caption{\small{The graph of $S(x)\equiv S_{2,2}(x)$ against $\log_{10}x$ for different ranges of $x$.}}
\end{center}
\end{figure}

The case of the original Riesz function has $m=p=2$ and from (\ref{e21}) we have
\[S_{2,2}(x)=\sum_{n=1}^\infty \frac{\mu(n)}{n^2} (1-e^{-X_n})-\frac{1}{\zeta(2)},\]
where $\zeta(2)=\pi^2/6$. The results are shown in Fig.~3 for the range $0\leq x\leq 10^8$ obtained using the scheme (\ref{e29}) with $p=2$. This function presents a similar behaviour with its value decreasing until about $10^{-8}$ before commencing to oscillate about the zero line. The first zero has the value $x\doteq 1.15671$. A plot of $x S_{2,2}(x)$ is also given in \cite{CW}. The successive maxima and minima in the oscillatory region are computed as the zeros of the derivative
\[S_{2,2}'(x)=\frac{1}{\zeta(4)}-\sum_{n=1}^\infty \frac{\mu(n)}{n^4} (1-e^{-X_n}).\]
The results of these calculations together with the corresponding values of $S(x_k^\pm)$ are shown in Table 2. Plots of $\log\,S_{2,2}(x_k^\pm)$ against $\log\,x$ are shown in Fig.~4, where it is seen that they reveal a linear variation to a good approximation. The dashed lines in these figures have slope equal to $-0.75$, thereby numerically confirming the estimate in (\ref{e12}).
\begin{table}[h]
\caption{\footnotesize{Values of the maxima and minima $\log_{10} x_k^\pm$ and the corresponding values of $S_{2,2}(x_k^\pm)$.}}
\begin{center}
\begin{tabular}{|c|ll|ll|}
\hline
&&&&\\[-0.3cm]
\mcol{1}{|c|}{$k$} & \mcol{1}{c}{$\log_{10} x_k^+$} & \mcol{1}{c|}{$S_{2,2}(x_k^+)$} &\mcol{1}{c}{$\log_{10} x_k^-$} & \mcol{1}{c|}{$S_{2,2}(x_k^-)$} \\
[.1cm]\hline
&&&&\\[-0.25cm]
1 & 4.48752 & $4.97204\times 10^{-8}$  & 4.31797 & $-6.46896\times 10^{-9}$ \\
2 & 4.87969 & $1.97351\times 10^{-8}$  & 4.69479 & $-1.67414\times 10^{-8}$ \\
3 & 5.26779 & $9.09065\times 10^{-9}$  & 5.07699 & $-1.09071\times 10^{-8}$ \\
4 & 5.65449 & $4.53355\times 10^{-9}$  & 5.46263 & $-5.99878\times 10^{-9}$ \\
5 & 6.04080 & $2.28418\times 10^{-9}$  & 5.84775 & $-3.14719\times 10^{-9}$ \\
6 & 6.42706 & $1.17580\times 10^{-9}$  & 6.23435 & $-1.62554\times 10^{-9}$ \\
7 & 6.81307 & $5.98676\times 10^{-10}$ & 6.61979 & $-8.37753\times 10^{-10}$ \\
8 & 7.19932 & $3.07392\times 10^{-10}$ & 7.00651 & $-4.29573\times 10^{-10}$ \\
\hline
\end{tabular}
\end{center}
\end{table}

\begin{figure}[th]
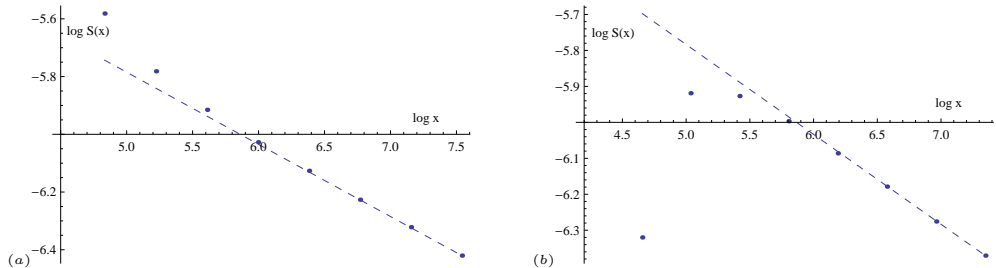

	\begin{center}	{\tiny($a$)}\includegraphics[width=0.4\textwidth]{RieszFig2a}
	\qquad
	{\tiny($b$)}\includegraphics[width=0.4\textwidth]{RieszFig2b} 
\caption{\small{Plots of (a) $\log_{10} S(x_k^+)$ and (b) $\log_{10} S(x_k^-)$  against $\log_{10}x$ (where $S(x)\equiv S_{2,2}(x)$). The dashed lines have slope $-0.75$.}}
\end{center}
\end{figure}

It has been shown in \cite{LBD} that $S_{2,2}(x)$ has an infinite number of zeros. It is very probable that $S_{2,1}(x)$ also has an infinite number of zeros. 

\vspace{0.6cm}

\begin{center}
{\bf 4.\ An asymptotic expansion}
\end{center}
\setcounter{section}{4}
\setcounter{equation}{0}
\renewcommand{\theequation}{\arabic{section}.\arabic{equation}}
An integral representation of $S_{m,p}(x)$ in the form of a Mellin-Barnes integral is given by
\bee\label{e41}
S_{m,p}(x)=-\frac{1}{2\pi i} \int_{c-\infty i}^{c+\infty i}\frac{\g(s) x^{-s}}{\zeta(p-ms)}\,ds\qquad (0<c<c_0),
\ee
where $c_0=(p-\fs)/m$. The integrand possesses simple poles at $s=0, -1 -2, \ldots$ and at the trivial zeros of the zeta function at $s=(p+2k)/m$, $k=1, 2, \ldots\,$. {\it On the assumption of the Riemann hypothesis}, there is also an infinite number of (simple) poles on the line $\Re (s)=c_0$ given by $s=c_0\pm i\gamma_k/m$ ($k=1, 2, \ldots$), where $\zeta(\fs\pm i\gamma_k)=0$. Assuming that it is permissible to displace the integration path past this line, we obtain the result
\[S_{m,p}(x)=-\frac{2 x^{-c_0}}{m} \Re \sum_{k=1}^\infty \frac{\g(c_0-i\gamma_k/m)}{\zeta'(\fs+i\gamma_k)}\ x^{i\gamma_k/m}+O(x^{-(p+2)/m})\]
as $x\to+\infty$. The details of the case $m=2$, $p=1$ are discussed in \cite[\S 2.5]{HL}; see also the account presented in \cite[p.~143]{PK}.

If we now set
\[A_k=\frac{\g(c_0-i\gamma_k/m)}{\zeta'(\fs+i\gamma_k)},\qquad \psi_k=\pi+\arg A_k,\]
we find that
\bee\label{e42}
S_{m,p}(x)=\frac{2x^{c_0}}{m} \sum_{k=1}^\infty |A_k| \cos\bl\{\frac{\gamma_k}{m}\,\log\,x+\psi_k\br\}+O(x^{-(p+2)/m}).
\ee
The convergence of the sum (\ref{e42}) is difficult to establish. The gamma function present in the coefficients $A_k$ decays very rapidly for increasing $k$, since from Stirling's formula it contains the exponential factor
$\exp\,[-\pi\gamma_k/(2m)]$ for large $\gamma_k$. The magnitude of $\zeta'(\fs+i\gamma_k)$ (which is non-zero on the assumption that the non-trivial zeros are all simple) generally increases with $k$, but it is possible that there are zeros for which this quantity could become small.

When $m=p=2$, we obtain the expansion
\bee\label{e43}
S_{2,2}(x)=x^{-3/4} \sum_{k=1}^\infty |A_k| \cos\,\bl\{\frac{\gamma_k}{2}\,\log\,x+\psi_k\br\} +O(x^{-2}) 
\ee
with $A_k=\g(\f{3}{4}-\fs i\gamma_k)/\zeta'(\fs+i\gamma_k)$.
The graph of $S_{2,2}(x)$ against $\log_{10}x$ compared with the expansion (\ref{e43}) truncated after $k=5$ terms is shown in Fig.~5. It is seen that for $x\,\gtwid\,10^6$ the curves are indistinguishable on the scale of the figure.
\begin{figure}[th]
	\begin{center}	\includegraphics[width=0.6\textwidth]{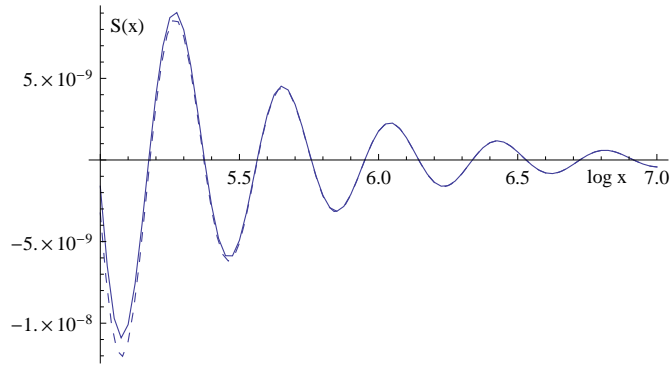}
\caption{\small{The graph of $S(x)\equiv S_{2,2}(x)$ against $\log_{10}x$ compared with the asymptotic form (\ref{e43}) (dashed curve).}}
\end{center}
\end{figure}
\vspace{0.6cm}

\begin{center}
{\bf 5.\ Concluding remarks}
\end{center}
\setcounter{section}{5}
\setcounter{equation}{0}
\renewcommand{\theequation}{\arabic{section}.\arabic{equation}}
The numerical results obtained in Section 3 are indicative only. It appears from the numerical investigation -- without any reference to the Riemann hypothesis -- that the large-$x$ behaviour of the functions $S_{2,1}(x)$ and $S_{2,2}(x)$ possesses respectively the $x^{-1/4}$ and $x^{-3/4}$ decay superimposed on an oscillatory structure.
A striking feature of the plots in Figs.~1 and 3 is the fact that the final decaying oscillatory structure is not obtained until $x$ has attained the value of approximately $10^5$. This is an unusual occurrence since most special functions begin to exhibit their asymptotic structure for often surprisingly modest values of the variable.

The computation of $S_{m,p}(x)$ for $m\geq 3$ is made easier since the cut-off value $N$ then scales like $x^{1/m}$ and the rate of decay of the various series in (\ref{e29}) is correspondingly more rapid. As an example, the case $m=3$, $p=2$ is shown in \cite{PRiesz}. The behaviour is found to be similar to that depicted in Figs,~1 and 3, with the maxima and minima following an approximate $x^{-1/2}$ scaling predicted by (\ref{e42}).
\vspace{0.6cm}

\begin{center}
{\bf Appendix: \ Derivation of Lemma 1}
\end{center}
\setcounter{section}{1}
\setcounter{equation}{0}
\renewcommand{\theequation}{\Alph{section}.\arabic{equation}}
From the contiguous relation satisfied by the confluent hypergeometric function ${}_1F_1(a;b,z)$ \cite[(13.3.2)]{DLMF}
\[b(b-1)  {}_1F_1(a;b-1;z)+b(1-b-z) {}_1F_1(a;b;z)+z(b-a) {}_1F_1(a;b+1;z)=0,\]
we obtain, with $a=1$, $b=k+2$, the recursion formula satisfied by $f_k(z):=e^{-z} {}_1F_1(1;k+1;z)$ in the form
\[f_k(z)-\bl(1+\frac{z}{k+1}\br) f_{k+1}(z)+\frac{z}{k+2} f_{k+2}(z)=0 \qquad (k\geq 0).\]
Repeated use of this result, combined with $f_0(z)=1$ and the partial sum of the exponential series $e_k(z)$ defined in (\ref{e24}), shows successively that
\begin{eqnarray*}
f_1(z)&=&\frac{1}{e_2(z)}+\frac{z e_1(z)}{2e_2(z)}\,f_2(z),\\
f_2(z)&=&\frac{1}{e_3(z)}+\frac{z e_2(z)}{3e_3(z)}\,f_3(z),
\end{eqnarray*}
and, in general,
\bee\label{a1}
f_{k-1}(z)=\frac{1}{e_k(z)}+\frac{z e_{k-1}(z)}{k e_k(z)}\,f_k(z) \qquad (k\geq 1).
\ee

Then, from (\ref{a1}) we find that
\begin{eqnarray*}
z f_1(z)&=&\frac{z}{e_2(z)}+\frac{z^2 f_2(z)}{2e_2(z)}=1-\frac{1}{e_2(z)}+\frac{z^2 f_2(z)}{2e_2(z)}\\
&=&1-\frac{1}{e_2(z)}+\frac{z^2}{2e_2(z)} \bl(\frac{1}{e_3(z)}+\frac{z e_2(z)}{3e_3(z)}\,f_3(z)\br)\\
&=&1-\frac{1}{e_3(z)}+\frac{z^3 f_3(z)}{3! e_3(z)},
\end{eqnarray*}
where we have used $e_1(z)=1$, $e_2(z)=1+z$ and have written $z^2/2!=e_3(z)-e_2(z)$. This procedure can be continued to produce the final result
\bee\label{a2}
zf_1(z)=1-\frac{1}{e_k(z)}+\frac{z^k f_k(z)}{k! e_k(z)}\qquad (k\geq 1).
\ee

\end{document}